\documentclass[12pt]{amsart}
\usepackage{amssymb}

\renewcommand{\a}{\alpha}
\renewcommand{\b}{\beta}
\renewcommand{\d}{\delta}

\newcommand{\e}{\varepsilon}
\newcommand{\g}{\gamma}

\newcommand{\var}{\varphi}
\newcommand{\s}{\sigma}

\newcommand{\cB}{\mathcal B}
\newcommand{\cC}{\mathcal C}
\newcommand{\cBst}{{\mathcal B}_{\pmb{s},\pmb{t}}}
\newcommand{\cCst}{{\mathcal C}_{\pmb{s},\pmb{t}}}
\newcommand{\ps}{\pmb{s}}
\newcommand{\pt}{\pmb{t}}
\newcommand{\pu}{\pmb{u}}
\newcommand{\pv}{\pmb{v}}
\newcommand{\pw}{\pmb{w}}
\newcommand{\px}{\pmb{x}}
\newcommand{\py}{\pmb{y}}
\newcommand{\pz}{\pmb{z}}
\newcommand{\Pst}{P_{\pmb{s},\pmb{t}}}
\newcommand{\Qst}{Q_{\pmb{s},\pmb{t}}}

\newcommand{\bell}{\mbox{\boldmath$\ell$}}

\newtheorem{theorem}{Theorem}[section]
\newtheorem{mtheorem}{Main Theorem}[section]
\newtheorem{lemma}{Lemma}[section]
\newtheorem{corollary}{Corollary}[section]
\newtheorem{proposition}{Proposition}[section]
\newtheorem{notation}{Notation}[section]
\newtheorem{definition}{Definition}[section]
\newtheorem{question}{Question}[section]

\begin{document}

\title[Hyperbolic group $C^*$-algebras]{Hyperbolic group
$C^*$-algebras and free-product $C^*$-algebras as compact
quantum metric spaces}
\author{Narutaka Ozawa and Marc A. Rieffel}
\address{Department of Mathematical Science, 
University of Tokyo, Komaba, 153-8914, Japan}
\email{narutaka@ms.u-tokyo.ac.jp}
\address{
Department of Mathematics \\
University of California \\ Berkeley, CA 94720-3840, U.S.A.}
\email{rieffel@math.berkeley.edu}
\date{April 17, 2003}
\thanks{The first author was supported by the Japan Society for 
the Promotion of Science Postdoctral Fellowships 
for Research Abroad, 
and the research of the second author was 
supported in part by National Science Foundation grants DMS99-70509
and DMS-0200591.}

\begin{abstract}
Let $\ell$ be a length function on a group $G$, and let $M_{\ell}$ denote the
operator of pointwise multiplication by $\ell$ on $\bell^2(G)$. 
Following Connes,
$M_{\ell}$ can be used as a ``Dirac'' operator for $C_r^*(G)$.  It defines a
Lipschitz seminorm on $C_r^*(G)$, which defines a metric on the state space of
$C_r^*(G)$. We show that if $G$ is a hyperbolic group and if $\ell$ is
a word-length function on $G$, then the topology from this metric 
coincides with the
weak-$*$ topology (our definition of a ``compact quantum metric 
space''). We show that a convenient framework is that of filtered
$C^*$-algebras which satisfy a suitable `` Haagerup-type'' condition. We
also use this
framework to prove an analogous fact for certain reduced
free products of $C^*$-algebras.
\end{abstract}

\maketitle
\allowdisplaybreaks

\setcounter{section}{-1}
\section{Introduction}
\label{sec0}

The group $C^*$-algebras of discrete groups provide a much-studied class of
``compact non-commutative spaces'' (that is, unital $C^*$-algebras). 
In \cite{Cn1} Connes showed that the ``Dirac'' operator of a spectral
triple (i.e. of an unbounded 
Fredholm module) over a unital $C^*$-algebra provides in a natural way 
a metric on the state space of the algebra.  The
class of examples most discussed in \cite{Cn1} consists of the  group 
$C^*$-algebras of discrete groups, with the Dirac operator coming in a 
simple way from a word-length function on the group. 
In \cite{Rf1}, \cite{Rf2}  
the second author pointed out that, motivated by what happens 
for ordinary compact metric spaces, it is natural to desire that 
for a spectral triple the topology from the metric on
the state space coincides with the weak-$*$ topology (for which the 
state space is compact).  This property was verified in \cite{Rf1} 
for certain examples. In \cite{Rf3} this property was taken as the
defining property for a ``compact quantum metric space''.

In \cite{Rf4} the second author studied this property for Connes'
original example of discrete groups with Dirac operators coming
from a word-length functions, but was able to verify this property
only for the case when the group is ${\mathbb Z}^n$. This already
took a long and interesting argument. We refer the reader to the
introduction of \cite{Rf4} for a more extensive discussion of this
whole matter.

In the present paper we verify the property for the case of
hyperbolic discrete groups. In the course of studying this case
we discovered that a natural setting was that of filtered $C^*$-algebras
with faithful trace. Voiculescu had shown earlier \cite{Vo2} how
to define an appropriate Dirac operator in that setting. In Section 1
we formulate in that setting a ``Haagerup-type condition'', which
in Sections 2 and 3 we show is sufficient to imply that the metric 
from the Dirac operator gives the state space the weak-$*$ topology.
Then in Section 4 we show that this Haagerup-type condition is
satisfied in the case of hyperbolic groups. We mention that quite
recently Antonescu and Christensen \cite{AC} showed that for 
non-Abelian free groups the metric on the state space gives the
state space finite diameter. Their techniques are close to ours,
but make explicit the relationship with Schur multipliers.

In Section 5 we show that the Haagerup-type condition fails for
the groups ${\mathbb Z}^n$ for $n \geq 2$ with their standard length
functions, and for groups which contain an amenable group of growth
$\geq 4$ for the length function in use. 
Since the approach used in the present
paper is entirely different from that used in \cite{Rf4} to 
successfully treat ${\mathbb Z}^n$, this raises the interesting question 
of finding a unified approach which covers both cases. And there
remains wide open the question of what happens for other classes
of groups, such as the discrete Heisenberg group and other nilpotent discrete
groups.

Finally, in Section 6 we show that the Haagerup-type condition is
satisfied by the reduced free product of any two filtered
$C^*$-algebras which satisfy the Haagerup-type condition. (Their
filtrations give in a natural way a filtration on the free product.)
This provides yet more examples of compact quantum metric spaces. 

We are very much indebted to Gilles Pisier for giving us a proof 
that for the free group on $n$ generators with its standard word-length 
function the corresponding metric on the state space gives the state 
space finite diameter. This showed us how to begin proving things in 
the direction which we have pursued here. We also warmly thank him 
for valuable comments on our manuscript.

\section{Filtered $C^*$-algebras}
\label{sec1}

We let $A$ be a unital $*$-algebra over ${\mathbb C}$ which has a
$*$-filtration
$\{A_n\}$ by finite-dimensional subspaces.  Just as in \cite{Vo2}
this means that
$A_m \subset A_n$ if $m < n$, $A = \cup_{n=0}^{\infty} A_n$, $A_n^*
= A_n$ and
$A_mA_n \subseteq A_{m+n}$, and $A_0 = {\mathbb C}1_A$.  We assume
further that we
are given a faithful state, $\s$, on $A$, that is, a linear
functional such that
$\s(a^*a) > 0$ for all $a \in A$ unless $a = 0$, and $\s(1_A) = 1$.
Let ${\mathcal
H} = L^2(A,\s)$ denote the corresponding GNS Hilbert space.  We
assume that the left
regular representation of $A$ on ${\mathcal H}$ is by bounded operators, 
and we
identify $A$ with the corresponding algebra of operators on
${\mathcal H}$.  We let
$\|\cdot\|$ denote the operator norm of $A$.  Our notation will not
distinguish between
$a$ as an operator on ${\mathcal H}$ and $a$ as a vector in
${\mathcal H}$, so the
context must be examined to see which is intended.  We let $\|a\|_2$
denote the norm
of $a$ as a vector in ${\mathcal H}$.

We can view each $A_n$ as a finite-dimensional, thus closed, subspace
of ${\mathcal
H}$.  We let $Q_n$ denote the orthogonal projection of ${\mathcal H}$
onto $A_n$.
We then set $P_n = Q_n - Q_{n-1}$ for $n \geq 1$, and $P_0 = Q_0$.   
The $P_n$'s are mutually
orthogonal, and $\sum P_n = I_{\mathcal H}$ for the strong operator
topology.  For
each $a \in A$ and each $n$ we set $a_n = P_n(a)$, where here $a$ is 
viewed as a
vector.  Then $a_n \in A_n$, but $a_n \notin A_{n-1}$ unless $a_n =
0$.  Furthermore
$a = \sum a_n$, with at most $p$ non-zero terms in the sum if $a \in A_p$.

For the above situation we define, as in \cite{Vo2}, an unbounded
operator, $D$, on
${\mathcal H}$ by $D = \sum_{n=1}^{\infty} nP_n$.  Notice that $A$ is
contained in
the domain of $D$.  The following lemma is part of proposition~$5.1$d
of \cite{Vo2}.
We include the proof here since we will need a similar argument in
Section~\ref{sec3}.

\begin{lemma}
\label{lem1.1}
For any $a \in A$ the operator $[D,a]$ has dense domain and is a
bounded operator.
\end{lemma}

\begin{proof}
Clearly $A$ is contained in the domain of $[D,a]$, and $A$ is dense.
Suppose that
$a \in A_p$.  Then for any given $m,n \ge 0$, if $P_maP_n \ne 0$ then
there is a
$\xi \in A_n$ such that $a\xi \in A_m$.  Since $A_pA_n \subseteq A_{p+n}$, it
follows that $p+n \ge m$.  On taking the adjoint, we see that
$P_ma^*P_n \ne 0$, so
that $p+m \ge n$.  Thus $|m-n| \le p$.  Consequently,
\[
a = \sum_{|m-n| \le p} P_maP_n,
\]
converging in the strong operator topology.  For each $j$ with $|j| \le p$ set
\[
T_j = \sum P_maP_{m-j}.
\]
Because the range of the terms $P_maP_{m-j}$ are orthogonal for fixed
$j$, as are
the ``domains'', we have
\[
\|T_j\| = \sup_m \|P_maP_{m-j}\| \le \|a\|.
\]
But for any $m,n \ge 0$ we have
\[
[D,P_maP_n] = (m-n)P_maP_n.
\]
In particular, $[D,P_maP_{m-j}] = jP_maP_n$.  Thus $[D,T_j] = jT_j$.
Since $a =
\sum T_j$, we obtain
\[
[D,a] = \sum_{|j| \le p} jT_j.
\]
\end{proof}

Thus $(A,{\mathcal H},D)$ is a spectral triple (or unbounded Fredholm
module) as
defined by Connes \cite{Cn1}, \cite{Cn2}.  We can then define a
seminorm, $L$, on
$A$ by
\[
L(a) = \|[D,a]\|.
\]
  From the proof of Lemma~\ref{lem1.1} we can see that $L$ will be a Lipschitz
seminorm on $A$ in the sense \cite{Rf2} that $L(a) = 0$ exactly if $a
\in {\mathbb
C}1_A = A_0$.

As pointed out by Connes, for any spectral triple $(A,{\mathcal
H},D)$, with $L$
defined as above, we can define a metric, $\rho_L$, on the state
space $S(A)$ of $A$
by
\[
\rho_L(\mu,\nu) = \sup\{|\mu(a)-\nu(a)|: L(a) \le 1\},
\]
(which may be $+\infty$).  As discussed in \cite{Rf1}, \cite{Rf2},
\cite{Rf3} it is
natural to ask whether the topology on $S(A)$ determined by $\rho_L$
agrees with the
weak-$*$ topology, as happens for ordinary compact metric spaces
$(X,\rho)$ and the
usual Lipschitz seminorm on $C(X)$.  If so, then \cite{Rf2} we call $L$ a
``Lip-norm''.  We consider a unital (pre-) $C^*$-algebra equipped
with a Lip-norm to
be a compact quantum metric space.

\setcounter{mtheorem}{1}
\begin{mtheorem}
\label{mthm1.1}
Let $A$, $\s$ and the $*$-filtration $\{A_n\}$ be as above, and let
$D$ and $L$ be
defined as above.  If furthermore there is a constant, $C$, such that
\[
\|P_ma_kP_n\| \le C\|a_k\|_2
\]
for all $a \in A$ and integers $m,n,k$, then $L$ is a Lip-norm.
\end{mtheorem}

As we will see at the end of Section 3, the key condition 
involving $C$ stated just above is closely
related to the Haagerup inequality. We will call a condition
of this kind a ``Haagerup-type condition''.

Necessary and sufficient conditions for a Lipschitz seminorm on a
pre-$C^*$-algebra
to be a Lip-norm are given in \cite{Rf1} (in a more general context).  For our
present purposes it is convenient to reformulate these conditions slightly.

\setcounter{proposition}{2}
\begin{proposition}
\label{prop1.2}
Let $L$ be a Lipschitz seminorm on a unital pre-$C^*$-algebra $A$,
and let $\s$ be a
state of $A$.  Then $L$ is a Lip-norm if and only if
\[
\{a \in A: L(a) \le 1 \mbox{ and } \s(a) = 0\}
\]
is a norm-totally-bounded subset of $A$.
\end{proposition}

\begin{proof}
We apply theorem $1.8$ of \cite{Rf1}.  Let $E = \{a\in A: L(a) \le 1
\mbox{ and }
\s(a) = 0\}$.  Suppose first that $E$ is totally bounded.  As in
theorem $1.8$ of
\cite{Rf1} let ${\mathcal L}_1 = \{a \in A: L(a) \le 1\}$, and let
${\tilde A} =
A/{\mathbb C}1_A$ with the quotient norm.  Let ${\tilde {\mathcal
L}}_1$ denote the
image of ${\mathcal L}_1$ in ${\tilde A}$.  For any $a \in
{\mathcal L}_1$
the element $a-\s(a)1_A$ is in $E$.  Thus the image of $E$ in
${\tilde A}$ coincides
with ${\tilde {\mathcal L}}_1$.  Thus if $E$ is totally bounded then
so is ${\tilde
{\mathcal L}}_1$.  But this is exactly the condition in theorem $1.8$
of \cite{Rf1}
for $L$ to be a Lip-norm.  Conversely, if $L$ is a Lip-norm so that ${\tilde
{\mathcal L}}_1$ is totally bounded, then a simple $2\e$-argument
shows that $E$ is
totally bounded.
\end{proof}

\section{The action of the one-parameter group}
\label{sec2}

In this section we consider a Hilbert space ${\mathcal L}$ with a
sequence $\{P_n\}$
of mutually orthogonal projections whose sum is $I_{\mathcal H}$,
much as above.  We
set $D = \sum nP_n$, and for each $t \in {\mathbb R}$ we let $U_t =
e^{itD} = \sum
e^{itn} P_n$.  We let $\a_t$ denote the inner automorphism of ${\mathcal
B}({\mathcal H})$ defined by $\a_t(T) = U_tTU_t^*$.  Because the
spectrum of $D$
consists of integers, we can view $\a$ as an action of the circle
group ${\mathbb T} = {\mathbb R}/(2\pi {\mathbb Z})$.  
In general the function $t \mapsto \a_t(T)$ will not be
norm-continuous.  But
it is always strong-operator continuous.  Thus for any finite measure $\mu$ on
${\mathbb T}$ and any $T \in {\mathcal B}({\mathcal H})$ we can
define $\a_{\mu}(T)$
by
\[
(\a_{\mu}(T))\xi = \int_{\mathbb T} \a_t(T)\xi\  d\mu(t)
\]
for each $\xi \in {\mathcal H}$.  Then $\|\a_{\mu}(T)\| \le \|T\|
\|\mu\|_1$, where
$\|\mu\|_1$ is the total-variation norm.  Notice then that for any
$m,n \ge 0$ we
have
\begin{eqnarray*}
P_m\a_{\mu}(T)P_n &= &\int e^{imt} P_mTP_ne^{-int}d\mu(t) \\
&= &{\hat \mu}(n-m)P_mTP_n,
\end{eqnarray*}
where ${\hat \mu}$ is the Fourier transform of $\mu$.  In particular,
if $[D,T]$ is
a bounded operator, then
\begin{eqnarray*}
P_m \a_{\mu}([D,T])P_n &= &{\hat \mu}(n-m)P_m[D,T]P_n \\
&= &(m-n){\hat \mu}(n-m)P_mTP_n.
\end{eqnarray*}

For any integer $N \ge 0$ let $\var_N \in \ell^2({\mathbb Z})$ be defined by
$\var_N(k) = -1/k$ if $|k| > N$ and $0$ otherwise.  Then the inverse Fourier
transform, $\check{\var}_N$, of $\var_N$ is in $L^2({\mathbb T})$, and so
in $L^1({\mathbb
T})$.  Thus as the measure $\mu$ above we can use $\check{\var}_N(t)dt$.
With some abuse
of notation we denote the corresponding operator by $\a_{\var_N}$.
For any $T \in
{\mathcal B}({\mathcal H})$ for which $[D,T]$ is bounded we set
\[
T^{(N)} = \a_{\var_N}([D,T]).
\]
Then for any $m,n \ge 0$ we have, as above,
\begin{eqnarray*}
P_mT^{(N)}P_n &= &(m-n)\var_N(n-m)P_mTP_n \\
&= &\left\{ \begin{array}{rl}
0 &\mbox{if $|m-n| \le N$} \\
P_mTP_n &\mbox{if $|m-n| > N$.}
\end{array} \right.
\end{eqnarray*}
Thus
\[
T^{(N)} = \sum_{|m-n| > N} P_mTP_n.
\]
Furthermore, 
\[
\|T^{(N)}\| \le 2\pi\|\var_N\|_2\|[D,T]\|,
\]
since 
$\|\check{\var}_N\|_1 \leq \sqrt{2\pi}\|\check{\var}_N\|_2
= 2\pi \|\var_N\|_2$.  Notice that
$\|\var_N\|_2
\rightarrow 0$ as $N \rightarrow +\infty$.

\section{The proof of the Main Theorem}
\label{sec3}

We resume the notation of Section~\ref{sec1}.  According to
Proposition~\ref{prop1.2} we must show that, under the hypotheses of the Main
Theorem, the set
\[
E = \{a \in A: \|[D,a]\| \le 1 \mbox{ and } \s(a) = 0\}
\]
is totally bounded in $A$ for the operator norm.  Given $a \in A$, we
set $a_n =
P_n(a)$ as in Section~\ref{sec1}, so that $a = \sum a_n$.  The
condition that $\s(a)
= 0$ is then just the condition that $a_0 = 0$.

Let $\e > 0$ be given.  We now show that $E$ can be covered by a
finite number of
$3\e$-balls.  For $\var_N$'s as in the previous section, choose $N$
large enough
that $2\pi\|\var_N\|_2 < \e$.  For $a \in E$ define $a^{(N)}$ as in the
previous section
by $a^{(N)} = \a_{\var_N}([D,a])$.  Then from the discussion there we have
$\|a^{(N)}\| < \e$.
Set $a^N = a - a^{(N)}$, so that $\|a-a^N\| < \e$. Since as above
\[
a^{(N)} = \sum_{|m-n| > N} P_maP_n,
\]
we have
\[
a^N = \sum_{|m-n| \le N} P_maP_n,
\]
which converges in the strong operator topology.  Note that in
general $a^N \notin
A$.

Let $1_A$ be viewed as a vector in $L^2(A,\s)$, so that $\|1_A\|_2 =
1$ and $D(1_A)
= 0$.  Then for any $a \in A$ we have
\[
[D,a](1_A) = D(a) = \sum na_n.
\]
Since the $a_n$'s are mutually orthogonal, it follows that for $a \in
E$ we have
\[
\sum n^2\|a_n\|_2^2 \le \|[D,a]\|^2 \le 1.
\]
Then from the Cauchy--Schwarz inequality we see that for any integer
$K \ge 0$ we
have
\begin{eqnarray*}
\sum_{n > K} \|a_n\|_2 &= &\sum_{n>K} (n^{-1})(n\|a_n\|_2) \\
&\le &\left( \sum_{n>K} n^{-2}\right)^{1/2}(\sum n^2\|a_n\|_2^2)^{1/2} \\
&\le &\left( \sum_{n>K} n^{-2}\right)^{1/2}.
\end{eqnarray*}

We now choose $K$ large enough that
\[
\left( \sum_{n>K} n^{-2}\right)^{1/2} < \e(C(2N+1))^{-1}.
\]
For each $a \in A$ set ${\hat a}_K = \sum_{k\le K} a_k$ and ${\tilde a}_K =
\sum_{k>K} a_k$, so that $a = {\hat a}_K + {\tilde a}_K$.  Then
\[
a = a^N + a^{(N)} = {\hat a}_K^N + {\tilde a}_K^N + a^{(N)},
\]
where $\hat a^N_K = (\hat a_K)^N$ and similarly for $\tilde a^N_K$.
For $a \in E$ we have chosen $N$ so that $\|a^{(N)}\| < \e$.  We show next that
$\{{\hat a}_K^N: a \in E\}$ is totally bounded.  Then we will show
that because of
our choice of $K$ we have $\|{\tilde a}_K^N\| < \e$ for any $a \in E$.  It will
follow immediately that $E$ can be covered by a finite number of
$3\e$-balls, as
desired.

For any $a \in E$ we have
\[
\|{\hat a}_K\|_2 \le \sum_{k \le K} \|a_k\|_2 \le \left( \sum_{k=1}^{\infty}
n^{-2}\right)^{1/2}.
\]
Thus $\{{\hat a}_K: a \in E\}$ is a bounded subset of the finite
dimensional vector
space $A_K$.  The map $a \mapsto {\hat a}^N$ is linear, and so when
restricted to
$A_K$ it must carry $\{{\hat a}_K: a \in E\}$ to a bounded subset of a
finite-dimensional subspace of ${\mathcal B}({\mathcal H})$.  Thus
$\{{\hat a}_K^N:
a \in E\}$ is totally bounded, as needed. (This is the only place
in this proof where we use the assumption that the $A_n$'s are
finite dimensional. Without that assumption this proof only shows 
that the metric on $S(A)$ gives $S(A)$ finite diameter.)

We now show that $\|{\tilde a}_K^N\| < \e$ for $a \in E$.  It is convenient to
first  show the following slightly more general fact:

\setcounter{lemma}{0}
\begin{lemma}
\label{lem3.1}
With notation as above, for any $a \in A$ we have
\[
\|a^N\| \le (2N+1)C \sum_{k=0}^{\infty} \|a_k\|_2.
\]
\end{lemma}

\begin{proof}
For each integer $j$ with $|j| \le N$ set,
\[
T_j = \sum_m P_maP_{m-j}.
\]
As in the proof of Lemma~\ref{lem1.1} we have
\[
\|T_j\| = \sup_m \|P_maP_{m-j}\|.
\]
For each integer $m$ we have, by hypothesis,
\[
\|P_maP_{m-j}\| \le \sum_k \|P_ma_kP_{m-j}\| \le C\sum \|a_k\|_2,
\]
so that $\|T_j\| \le C\sum \|a_k\|_2$.  Since $a^N = \sum_{|m-n| \le
N} P_maP_n =
\sum_{|j| \le N} T_j$, we obtain the asserted fact.
\end{proof}

Now for any $a \in E$, because ${\tilde a}_K = \sum_{k > K} a_k$, the above
proposition gives
\begin{eqnarray*}
\|{\tilde a}_K^N\| &\le &(2N+1)C \sum_{k>K} \|a_k\|_2 \\
&\le &(2N+1)C\left( \sum_{k>K} (k^{-2})\right)^{1/2} < \e
\end{eqnarray*}
by our choice of $K$, as needed. This concludes the proof of Main 
Theorem \ref{mthm1.1}.

We show next that from our Haagerup-type condition we can obtain a Haagerup
inequality in its more usual form.  Let $a \in A$, and let the
$a_k$'s be its components as above.  For any $k$ and for $|j| \le k$ set $T_j =
\sum P_ma_kP_{m-j}$, much as above.  Then, as above,
\[
\|T_j\| = \sup_m \|P_ma_kP_{m-j}\| \le C\|a_k\|_2.
\]
Since, as above, $a_k = \sum _{|j| \le k} T_j$, we obtain the following 
analog of the third line of the proof of lemma 1.4 of \cite{Hg},
which we record for later use:

\setcounter{lemma}{1}
\begin{lemma}
\label{lem3.2}
With notation as above, we have
\[
\|a_k\| \le C(2k+1)\|a_k\|_2.
\]
\end{lemma}

Then from the Cauchy--Schwarz inequality we obtain
\begin{eqnarray*}
\|a\| &\le &\sum_k \|a_k\| \le \sum_k C(2k+1)\|a_k\|_2 \\
&= &C \sum_k (1/(k+1))(k+1)(2k+1)\|a_k\|_2 \\
&\le &C\left( \sum_{p \ge 1} 1/p^2\right)^{1/2}
\left(\sum(2k^2+3k+1)^2\|a_k\|_2^2\right)^{1/2}.
\end{eqnarray*}
If we note that $2k^2 + 3k + 1 \le 2(k+1)^2$ for $k \ge 0$, 
and set $C' = 2C\left(\sum_{p \ge 1} 1/p^2\right)^{1/2}$,  we obtain
the following inequality, which is similar to the usual form \cite{Cn2} 
for the Haagerup inequality for groups:

\setcounter{proposition}{2}
\begin{proposition}
\label{prop3.3}
For any $a \in A$ we have
\[
\|a\| \le C' \big(\sum (1+k)^4\|a_k\|_2^2\big)^{1/2}.
\]
\end{proposition}

We now obtain a related inequality which we will need shortly.

\setcounter{proposition}{3}
\begin{proposition}
\label{prop3.4}
There is a constant, $C''$,  such that for any integer $p$ and any 
$a \in A_p$ we have
\[
\|a\| \leq C'' (p+1)^{3/2} \|a\|_2 .
\]
\end{proposition}

\begin{proof}
We use Lemma \ref{lem3.2} to calculate that
\begin{eqnarray*}
\|a\| & \leq & \sum^p_0 \|a_k\| \leq C\big(\sum^p_0 (2k+1)\|a_k\|_2\big) \\
& \leq & C\big(\sum^p_0 (2k+1)^2\big)^{1/2}\big(\sum^p_0\|a_k\|^2_2\big)^{1/2} 
\\
& \leq & 2C\big(\sum^p_0 (k+1)^2\big)^{1/2}\|a\|_2.
\end{eqnarray*}
But
\[
\sum^p_0 (k+1)^2 \leq \int^{p+1}_0 (t+1)^2 \ dt = (1/3)\big((p+2)^3 - 1\big).
\]
Absorbing several factors into the constant, we obtain the desired
inequality.
\end{proof}

\section{Hyperbolic groups}
\label{sec4}

In this section we show that our Main Theorem applies to word-hyperbolic
groups.  There are several equivalent definitions of what it means for a metric
space to be hyperbolic \cite{GhH}.  We will find the following version
well-suited to our purposes.

\setcounter{definition}{0}
\begin{definition}
\label{def4.1}
A metric space $(X,\rho)$ is {\em hyperbolic} if there is a constant $\d \ge 0$
such that for any four points $x,y,z,w \in X$ we have
\[
\rho(x,y) + \rho(z,w) \le \max\{\rho(x,z) + \rho(y,w),\rho(x,w) + \rho(y,z)\} +
\d.
\]
If it is important to specify $\d$, we say that $X$ is $\d$-hyperbolic.
\end{definition}

Let $G$ be a finitely generated discrete group, and let $S$ be a finite
generating subset for $G$, with $S = S^{-1}$.  Let $\ell$ be the word-length
function on $G$ determined by $S$, and let $\rho$ be the corresponding
left-invariant metric on $G$ defined by $\rho(x,y) = \ell(x^{-1}y)$.  Then $G$
is said to be hyperbolic if the metric space $(G,\rho)$ is hyperbolic.  It is
not difficult to show \cite{GhH} that this is independent of the choice of the
finite generating set $S$.

For any discrete group $G$ and any integer-valued length function $\ell$ on $G$ 
we obtain a
$*$-filtration $\{A_n\}$ of the convolution algebra $A = C_c(G)$ of
complex-valued functions of finite support on $G$ by setting
\[
A_n = \{f \in A: f(x) = 0 \mbox{ if } \ell(x) > n\}.
\]
The involution on $A$ is defined, as usual, by $f^*(x) = (f(x^{-1}))^-$.  We
define a faithful trace, $\s$, on $A$ by $\s(f) = f(e)$, where $e$ denotes the
identity element of $G$.  The resulting GNS Hilbert space is $\bell^2(G)$, and
the left regular representation of $A$ on $\bell^2(G)$ is by bounded 
operators. The $C^*$-algebra generated by the left regular representation
is the reduced $C^*$-algebra of $G$, $C^*_r(G)$.
Thus we are in the setting of Section~\ref{sec1}. (With a bit of care with 
the bookkeeping, all the above applies also to the convolution algebra 
of $G$ twisted by 
a 2-cocycle, in the way that was explicitly carried out in \cite{Rf4}. Our 
results below also work for this case too.) 

The Dirac operator corresponding to the filtration is just
the operator $M_{\ell}$ of pointwise multiplication by $\ell$ on $\bell^2(G)$.
We can then define the seminorm $L$ on $A$ by $L(f) = \|[D,f]\|$, where $f$ on
the right is viewed as the convolution operator on $\ell^2(G)$.  We 
can then ask
whether $L$ is a Lip-norm.  Our Main Theorem provides a possible
tool for giving an affirmative
answer to this question. 

\setcounter{definition}{1}
\begin{definition}
\label{def4.2}
Let $\ell$ be an integer-valued length function on a group $G$. We say that 
$(G, \ell)$
satisfies a Haagerup-type condition if, for the filtration of
$C_c(G) \subseteq C^*_r(G)$ defined above, with its canonical trace,
the main condition of Theorem \ref{mthm1.1} is satisfied.
\end{definition}

\setcounter{proposition}{2}
\begin{proposition}
\label{prop4.3}
Let $G$ be a word-hyperbolic group, and let $\ell$ be the 
word-length function for a
finite generating subset of $G$. Then $(G, \ell)$ satisfies a
Haagerup-type condition.
\end{proposition}

\begin{proof}
A proof is 
essentially contained within Connes' proof of the Haagerup
inequality for hyperbolic groups given on page~241 of \cite{Cn2}.  
But since some significant details
are not included there, we give a complete proof here. The special case
of this proposition for the free group on finitely many
generators with its standard
word-length function relative to the given generators is explicitly given
by Haagerup as lemma 1.3 in \cite{Hg}, with $C = 1$. (See also lemma 1.1 of 
\cite{FT}, where it is remarked right after the proof of theorem 1.3 
that it also works
for the free group with countably many generators. But with an infinite
number of generators the subspaces $A_n$ of the filtration are
infinite dimensional, and so the proof of our Main Theorem 
\ref{mthm1.1} only shows that the state space has finite diameter.)

For any integer $j \ge 0$ let $E_j = \{x \in G: \ell(x) = j\}$.  
We must find a
constant, $C$, such that for any integers $k,m,n$, and any $f$ supported on
$E_k$ we have $\|P_mfP_n\| \le C\|f\|_2$.  This means that for any $\xi$
supported on $E_n$ we must have
\[
\left( \sum_{x \in E_m} |(f*\xi)(x)|^2\right)^{1/2} \le C\|f\|_2\|\xi\|_2.
\]
We examine $(f*\xi)(x)$.  Let $\d$ be a constant for which $G$, 
equipped with the
metric from $\ell$, is $\d$-hyperbolic as in Definition~\ref{def4.1}.  Now
\[
(f*\xi)(x) = \sum_{yz = x} f(y)\xi(z).
\]
If $(f*\xi)(x) \ne 0$ there must be some $y,z \in G$ such that
$x = yz$ with $\ell(y) = k$, $\ell(z) = n$,
and so if $x \in E_m$ we must have $m \le k+n$.  But also $z = y^{-1}x$, so we must
have $n \le k+m$, and so $|m-n| \le k$.  In the same way we obtain $|n-k| \leq m$. 
Let $p = k+n-m$.  If $p$ is even set
$q = p/2$, while if $p$ is odd set $q = (p-1)/2$.  In either case set ${\tilde
q} = p-q$, and notice that $q \le {\tilde q} \le q+1$.  Then $m =
(k-q)+(n-{\tilde q})$, and from $|m-n| \leq k$ 
it is easy to check that $k-q \ge 0$, while from $|n-k| \leq m$ it is
easy to check that $n-{\tilde
q} \ge 0$.  Consequently, for each $x \in E_m$ we can
choose ${\bar x},{\tilde x} \in G$ such that $x = {\bar x}{\tilde x}$ and
$\ell({\bar x}) = k-q$, while $\ell({\tilde x}) = n-{\tilde q}$. 
This choice is
usually not unique, but we fix it for the rest of the proof.

Suppose now that $x \in E_m$ and $x=yz$ for some $y \in E_k$ and $z \in E_n$.
We apply Definition \ref{def4.1} to the four points $(e,x,{\bar 
x},y)$ to obtain
\[
\rho(e,x) + \rho(y,{\bar x}) \le \max\{\rho(e,{\bar x}) + \rho(y,x), \ 
\rho(e,y) + \rho(x,{\bar x})\} + \d.
\]
But $\rho(e,{\bar x}) + \rho(y,x) = (k-q)+n$, while $\rho(e,y) + \rho(x,{\bar
x}) = k+(n-{\tilde q})$.  Consequently
\[
\rho(y,{\bar x}) \le k-q+n-m+\d = {\tilde q} + \d.
\]
Thus $y = {\bar x}u$ for some $u$ 
with $\ell(u) \le {\tilde q} + \d$.  Then $z = y^{-1}x = 
u^{-1}{\bar x}^{-1}x = u^{-1}{\tilde x}$.  Since this is true for
all such $x,y$, we see that
\[
(f*\xi)(x) = \sum \{f({\bar x}u)\xi(u^{-1}{\tilde x}): \ell(u) \le {\tilde q} +
\d\}.
\]
We can apply the Cauchy--Schwarz inequality to this to get
\begin{eqnarray*}
&   & |(f*\xi)(x)|^2  \\
&\le &\left(\sum \{|f({\bar x}u)|^2: \ell(u) \le {\tilde q} +
\d\}\right)\left(\sum \{|\xi(v{\bar x})|^2: \ell(v) \le {\tilde q}+\d\right).
\end{eqnarray*}

For any $y \in E_k$ let us consider how many decompositions there are of the
form $y = su$ such that $\ell(s) = k-q = \ell({\bar x})$ and $\ell(u) \le
{\tilde q} + \d$.  Let $y=tw$ be another such decomposition.  We apply
Definition \ref{def4.1} to the four points $e,y,s,t$ to obtain
\[
\rho(e,y) + \rho(s,t) \le \max\{\rho(e,s) + \rho(t,y), \ \rho(e,t) + 
\rho(s,y)\} + \d.
\]
But $\rho(e,s) + \rho(y, t) = k-q + \tilde q + \d = \rho(e,t) + \rho(s,y)$. 
It follows that
$k+\rho(s,t) \le {\tilde q} - q + k + 2\d$, so that $\rho(s,t) \le 1 + 2\d$.
In the same way we find that for any two factorizations $z = vs = wt$ with
$\ell(s) = \ell(t) = \ell({\tilde x}) = n-{\tilde q}$ and $\ell(v),\ell(w) \le
{\tilde q} + \d$ we have $\rho(v,w) \le 2\d$.  

Let $C$ be the number of elements of
$G$ in a ball of radius $1+2\d$.  Then 
the number of different $s$'s which can enter as above into the factorization
of $y$ is no larger than $C$, and thus the number of $u$'s is also no
larger than $C$. Similarly, the number of $v$'s which can enter as above 
into the factorization of $z$ is no larger than $C$.

We now claim that $\|f*\xi\|_2 \le C\|f\|_2\|\xi\|_2$.  From our earlier
calculations we know that
\begin{eqnarray*}
\|f*\xi\|_2^2 &= &\sum_x |(f*\xi)(x)|^2 \\
&\le &\sum_x \left( \sum_{\ell(u) \le  
{\tilde q} + \d} |f({\bar x}u)|^2\right) \left( \sum_{\ell(v) \le 
{\tilde q}+\d}
|\xi(v{\tilde x})|^2\right),
\end{eqnarray*}
while of course
\[
(\|f\|_2\|\xi\|_2)^2 = \sum_{\substack{\ell(y)=k \\ \ell(z)=n}}
|f(y)|^2|\xi(z)|^2.
\]
Thus to obtain our desired inequality it suffices to show that for any pair
$(y,z)$ with $\ell(y) = k$ and $\ell(z) = n$ the number of $x$'s for 
which there
are a $u$ and $v$ with $\ell(u) \le {\tilde q} + \d$ and $\ell(v) \le 
{\tilde q}
+ \d$ such that $y = {\bar x}u$ and $z = v{\tilde x}$ is no greater than $C^2$.
But suppose we have such $x,u,v$.  Then $x = {\bar x}{\tilde x} =
yu^{-1}v^{-1}z$.  Given our earlier bound on the number of such $u$'s 
and $v$'s,
it is now clear that the number of such $x$'s is indeed bounded by $C^2$.
\end{proof}

\setcounter{corollary}{3}
\begin{corollary}
\label{cor4.4}
Let $G$ be a word-hyperbolic group, and let $\ell$ be the 
word-length function for a
finite generating subset of $G$. Then the metric on $S(C^*_r(G))$
coming from using $\ell$ as a Dirac operator gives $S(C^*_r(G))$
the weak-* topology.
\end{corollary}

\section{Failure of the Haagerup-type condition}
\label{sec5}

In this section we show that the Haagerup-type condition 
often fails for groups which contain a copy of ${\mathbb Z}^d$ 
for $d \geq 2$, or other amenable groups with suitable growth. 
We begin with the following observation.

\setcounter{proposition}{0}
\begin{proposition}
\label{prop5.1}
Let $\ell$ be a length function on a group $G$, and let $\ell_H$ denote
the restriction of $\ell$ to a subgroup $H$. If $(G, \ell)$ satisfies
a Haagerup-type condition, then so does $(H,\ell_H)$. 
\end{proposition}

\begin{proof}
Since $G$ is the disjoint union of right cosets of $H$, the restriction
to $H$ of the left regular representation of $G$ is a direct sum of
copies of the left regular representation of $H$. Thus $C^*_r(H)$ is
isometrically embedded in $C^*_r(G)$. The restriction to $C^*_r(H)$ of the
canonical trace on $C^*_r(G)$ is the canonical trace on $C^*_r(H)$. The 
filtration of $C^*_r(H)$ for $\ell_H$ is just the intersection of
$C^*_r(H)$ with the filtration of $C^*_r(G)$ for $\ell$. The desired
conclusion follows easily.  
\end{proof}

\setcounter{proposition}{1}
\begin{proposition}
\label{prop5.2}
The group ${\mathbb Z}^2$ with the word-length function for its
standard basis does not satisfy a Haagerup-type condition. Thus
neither does ${\mathbb Z}^d$ for $d > 2$ with its standard word-length
function.
\end{proposition}

\begin{proof}

For ${\mathbb Z}^2$ and the standard word-length function $\ell$, given by
$\ell((p,q)) = |p| + |q|$, we need to show that there is no constant $C$ such
that $\|P_mfP_n\| \le C\|f\|_2$ for all $m,k,n$, where $f$ is supported on
$E_k$.  Let $k > 0$ be fixed, choose $n > k$, and set $m=n+k$.  Let $f$ be the
function which has value $(1/k)$ on the points $(p,k-p)$ of $E_k$ for which $1
\le p \le k$, and value $0$ elsewhere.  In the evident way we will consider $f$
to be a function just of $p$ when convenient.  Notice that $\|f\|_1 = 1$, so
that $\|P_mfP_n\| \le 1$, while $\|f\|_2 = 1/\sqrt{k}$.  Similarly, 
let $\xi$ be
the function which has value $1/\sqrt{n}$ on the points $(q,n-q)$ of $E_n$ for
which $1 \le q \le n$, and value $0$ elsewhere.  We can consider $\xi$ as a
function just of $q$.  Note that $\|\xi\|_2 = 1$.  We estimate
$\|P_mfP_n\xi\|$.  We will evaluate only on the points $(r,m-r)$ of $E_m$ for
which $k \le r \le n$.  Then with this restriction,
\begin{eqnarray*}
(P_mfP_n\xi)(r,m-r) &= &\sum_{1 \le p \le k} f(p)\xi(r-p) \\
&= &k(1/k)(1/\sqrt{n}) = 1/\sqrt{n}.
\end{eqnarray*}
Thus $\|P_mfP_n\xi\|_2^2 \ge (n-k)/n$, so that $\|P_mfP_n\| \ge
((n-k)/n)^{1/2}$.  Notice that this approaches $1$ as $n \rightarrow +\infty$.
But we could have chosen $k$ as large as desired, so that $\|f\|_2 = 
1/\sqrt{k}$
is as small as desired.  Thus there is no constant $C$ such that 
$\|P_nfP_m\| \le
C\|f\|_2$ for all $m,k,n$, where $f$ is supported on $E_k$.
\end{proof}

This, of course, raises the question of whether there is a way to 
give a unified
proof of both the Corollary \ref{cor4.4}  for hyperbolic groups and the
corresponding result in \cite{Rf4} for ${\mathbb Z}^d$, as well as the question
of what happens for other groups. Perhaps the ``bolic'' groups of
Kasparov and Skandalis \cite{KS} \cite{BK} provide a good class of
groups for which one might hope to find a unified proof.

Suppose now that $G$ is an amenable group, so that $C^*_r(G) = C^*(G)$.
Then the trivial representation of $G$ gives a representation of
$C^*_r(G)$. By using the trivial representation we see that if
$f \in C_c(G)$ and if $f \geq 0$ as a function, then $\|f\| = \|f\|_1$.
For each integer $p$ let $B_p = \{x \in G: \ell(x) \leq p\}$, and let
$\chi_p$ denote the characteristic function of $B_p$. Suppose that
$G$ satisfies a Haagerup-type condition. Then according to Proposition
\ref{prop3.4} there is a constant, $C'$, such that
\[
\| \chi_p\|_1 = \| \chi_p\| \leq C'(p+1)^{3/2} \| \chi_p\|_2.
\]
Let $|B_p|$ denote the number of elements in $B_p$. Then it follows
that $|B_p| \leq C'(p+1)^{3/2}|B_p|^{1/2}$. From this we obtain:

\setcounter{proposition}{2}
\begin{proposition}
\label{prop5.3}
Let $G$ be an amenable group, and let $\ell$ be an integer-valued 
length-function on $G$. If $(G,\ell)$ satisfies a Haagerup-type 
condition, then there is a constant, $C'$, such that for every
$p$ we have
\[
|B_p| \leq C'(p+1)^3.
\]
\end{proposition}

We now recall some well-known definitions and facts. (See page 12 of
\cite{GhH}.) For an integer-valued
length-function on $G$ we say that its rate of growth is polynomial
if there is an integer $n$ and a constant $C$ 
such that $|B_p| \leq Cp^n$ for all large enough $p$. We call the
smallest such $n$ the ``growth rate'' of $G$ for $\ell$. If $|B_p|$ grows
at a faster than polynomial rate, then we say that the growth rate of 
$G$ for $\ell$ is $\infty$. 

The idea of comparing the 2-norm with the 1-norm 
came from \cite{Jol}, where Jolissaint showed 
that an amenable group with the property (RD) is 
of polynomial growth. 

Let $S$ be a finite generating set for
$G$, and let $\ell_S$ be the corresponding word-length function. For
any length function $\ell$ on $G$ set $M = \max\{\ell(s): s \in S\}$.
Then it is easily seen that $\ell \leq M\ell_S$. Consequently the
growth rate of $G$ for $\ell$ is no smaller than that for $\ell_S$. In
particular, the growth rates of $G$ for any two word-length functions 
coincide. This common growth rate is called the growth rate of a given
finitely generated group. From the above observations 
and Proposition \ref{prop5.3} we obtain:

\setcounter{corollary}{3}
\begin{corollary}
\label{lem5.4}
If $G$ is a finitely generated amenable group, and if $G$ satisfies
a Haagerup-type condition for some length function, then the growth
rate of $G$ is no greater than 3.
\end{corollary}

\setcounter{corollary}{4}
\begin{corollary}
\label{cor5.5}
Let $G$ be any discrete group. If $G$ contains a finitely generated
amenable group whose growth rate is $\geq 4$, then there does not
exist a length function $\ell$ on $G$ such that $(G,\ell)$ satisfies
a Haagerup-type condition.
\end{corollary}

\setcounter{corollary}{5}
\begin{corollary}
\label{cor5.6}
If a group $G$ contains either ${\mathbb Z}^4$ or the discrete
Heisenberg group, then there does not
exist a length function $\ell$ on $G$ such that $(G,\ell)$ satisfies
a Haagerup-type condition.
\end{corollary}

\begin{proof}
Both ${\mathbb Z}^4$ and the discrete Heisenberg group have a growth
of 4. (See section 18 of chapter 1 of \cite{GhH} for the proof of this
for the Heisenberg group.)
\end{proof}

\setcounter{question}{6}
\begin{question}
\label{lem5.7}
Suppose that a group $G$ admits a finite generating set for
whose word-length function $\ell$ the pair $(G,\ell)$
satisfies a Haagerup-type condition. Must the group then 
be hyperbolic?
\end{question}

\section{Free-product $C^*$-algebras}
\label{sec6}

In this section we show that Main Theorem~\ref{mthm1.1}
applies to certain reduced free-product $C^*$-algebras.
Jolissaint \cite{Jol} showed that the property (RD) 
is preserved under forming free products, 
but his proof apparently does not work in our situation. 
Thus, we need a finer classification of types of words,
which unfortunately complicates the notation.

Let $A^1$ and $A^2$ be unital pre-$C^*$-algebras with
filtrations $\{A^1_m\}$ and $\{A^2_m\}$ respectively.
Let $A=A^1*A^2$ be the algebraic free product, with its
evident involution. We define a filtration (respecting
the involution) on $A$ by setting $A_n$ to be the linear
span of all products $A^{i_1}_{n_1}\cdots A^{i_\a}_{n_\a}$
with each $i_j = 1, 2$, with $i_j \ne i_{j+1}$ for 
$1 \leq j \leq \a - 1$, and with $\sum n_j \leq n$.

Let $\s^1$ and $\s^2$ be faithful tracial states on
$A^1$ and $A^2$. We let $\s = \s^1 * \s^2$ be the
corresponding faithful tracial state on $A$ which is
used to define \cite{Vo1} \cite{VDN} \cite{Av} the 
reduced free-product $C^*$-algebra structure on $A$.
Its defining properties are that its restrictions 
to $A^1$ and $A^2$ coincide with $\s^1$ and $\s^2$,
and that $\s(a^{i_1}_1 \cdots a^{i_\a}_\a) = 0$ if
$\s^{i_j}(a^{i_j}_j) = 0$ for all $j = 1, \ldots, \a$
and $i_j \neq i_{j+1}$ for $j = 1, \ldots, \a-1$.
The reduced $C^*$-norm on $A$ (for $\s_1$ and $\s_2$)
is then the operator norm for the GNS representation 
for $\s$ on $L^2(A, \s)$.

\setcounter{theorem}{0}
\begin{theorem}
\label{thm6.1} 

If $(A^1, \s^1)$ and $(A^2, \s^2)$ both satisfy a
Haagerup-type condition with constant $C$, then
$(A^1*A^2, \s^1*\s^2)$ satisfies a Haagerup-type
condition with constant $\sqrt{5}C$.
\end{theorem}

We remark that there are many examples to which this
theorem applies. In addition to the reduced group
$C^*$-algebras of hyperbolic groups studied in the
earlier sections of this paper, one can take any
finite-dimensional $C^*$-algebras with any filtrations.

This theorem is related to lemma 3.3 of \cite{DHR},
but in \cite{DHR} the algebras $A^1$ and $A^2$
are not assumed to be filtered, and so our
situation is substantially different from that
considered there.

We now establish some notation which will be used in
the proof. As in Section 1 we let $\{P^i_n\}$ be the 
family of mutually orthogonal projections corresponding
to the filtration $\{A^i_n\}$, for $i = 1, 2$, and we
let $\{P_n\}$ be the corresponding family on $A$ for
$\{A_n\}$. We let $E^i_n$ denote the range of $P^i_n$,
and similarly for $E_n$. Thus $E_0$ is the span of $1$,
while if $n \geq 1$ then $E_n$ is the orthogonal sum
of the spans of products 
$E^{i_1}_{n_1}\cdots E^{i_\a}_{n_\a}$ such that 
$n_j \geq 1$ for all $j$ and $i_j \neq i_{j+1}$ for
$j = 1, \ldots, \a - 1$ while $\sum n_j = n$. In order
to reduce notational clutter we will often omit the 
superscripts when they can be inferred from the context.
In particular, we will let $P_0^\perp$ denote the 
projection onto the orthogonal complement of $1$
for all three algebras.

Much as in section 2 of \cite{DHR} we choose for
$i = 1, 2$ an orthonormal basis $\cB^i_n$
for each $E^i_n$, with $\{1\}$ as the basis for 
$E^i_0$. But for convenience we also require that each 
basis element be self-adjoint. We can do this because
$\s_i$ is tracial. We let 
$\cB^i = \bigcup_n {\mathcal B}^i_n$, so that
$\cB ^i$ is a basis for $A^i$. We define $\ell$ on
each $\cB^i$ by $\ell(x) = n$ if $x \in \cB^i_n$.
For $x \in (\cB^1 \cup \cB^2)$
we define $\mu$ by $\mu(x) = i$ if $x \in \cB^i$,
and we define $\nu$ by $\nu(x) = i$ if $x \notin \cB^i$.
As in \cite{DHR} we obtain from $\cB^1$ and
$\cB^2$ an orthonormal basis $\cB$ for $A$. An element
of $\cB$ will be either $1$, or a product 
$\px = x_1 \cdots x_\a$ with 
$x_i \in (B^1 \cup B^2) \setminus \{1\}$ for each $i$
while $\mu(x_i) \neq \mu(x_{i+1})$ for 
$i = 1, \ldots, \a - 1$. We extend the definitions of
$\mu$ and $\nu$ to $\cB \setminus \{1\}$ by setting
$\mu(\px) = \mu(x_1)$ and $\nu(\px) = \nu(x_1)$ for
any $\px \neq 1$. Although $\mu(1)$ is undefined
(because $A$ is really the free product amalgamated
over ${\mathbb C}1$), we
will make the unusual convention that both 
$\mu(\px) = \mu(\py)$ and $\mu(\px) \neq \mu(\py)$
are simultaneously true if $\py = 1$. We set
$\ell(\px) = \sum \ell(x_j)$, with $\ell(1) = 0$.
We then set $\cB_n = \{\px: \ell(\px) = n\}$, and note
that $\cB_n$ is an orthonormal basis for $E_n$. (But
we note also that the elements of $\cB_n$ need not
be self-adjoint, though the involution carries $\cB_n$
into itself.) We will often write an
element $a$ of $E_n$ as 
$a = \sum_{\px \in \cB_n} a(\px)\px$.

Our objective is to show that for any $a \in E_k$
and any $m, n$ we have 
$\| P_m a P_n\|_2 \leq \sqrt{5}C\|a\|_2$, where on
the left side $a$ is viewed as an operator on
$L^2(A, \s)$. Thus we must show that if $\xi \in E_n$ 
then
$$
\|P_m(a\xi)\|_2 \leq \sqrt{5}C\|a\|_2 \|\xi\|_2.
$$
So we now fix $m, k$, and $n$ for the rest of the proof.
We can assume that $m, k$ and $n$ are all $\geq 1$,
since the desired inequality is very easily verified if 
any one of them is $0$. Somewhat as in Section 4 we
set $q = (k + n - m)/2$, but now $q$ need not be
an integer. Some of the objects considered below will
depend on $m, k$ and $n$, but to avoid notational clutter
we often will not indicate that dependence explicitly.

For any $a \in E_k$ we have 
$a = \sum_{\py \in \cB_k} a(\py)\py$. In the same
way, for $\xi \in E_n$ we have 
$\xi = \sum_{\pz \in \cB_n} \xi(\pz)\pz$. We find it
notationally convenient to work with $a^*\xi$
instead of $a\xi$. Then
$$
a^*\xi = \sum _{\py, \pz} {\bar a}(\py)\xi(\pz)\py^* \pz.
$$
Thus we need information about $P_m(\py^* \pz)$. So we
need to see how $\py^* \pz$ can be expressed in terms
of ``reduced words''. Let $\py = y_1 \cdots y_\b$
and $\pz = z_1 \cdots z_\g$. If $\mu(\py) \neq \mu(\pz)$,
then $\py^* \pz$ is already a reduced word, and
$\py^* \pz \in \cB_{k+n}$. Otherwise, if 
$\mu(\py) = \mu(\pz) $
then there is some integer $\d \geq 1$ such that $y_i = z_i$
for $i < \d$ while $y_\d \neq z_\d$ (with the latter
including the possibility that $y_\d$ or $z_\d$ is
not present, i.e. $\b < \d$ or $\g < \d$). If $\d = 1$
then $y_1 \neq z_1$ so that $P_0(y_1z_1) = 0$, and
$$
\py^*\pz = y_\b \cdots y_2 P_0^\perp(y_1z_1)z_2 \cdots z_\g,
$$
which is a reduced word. If $\d >1$ then $P_0(y_iz_i) = 1$
for $i < \d$, and so
\begin{eqnarray}
\py^*\pz &=& y_\b \cdots y_2 P_0(y_1z_1) z_2 \cdots z_\g
\ + \ y_\b \cdots y_2 P_0^\perp(y_1z_1) z_2 \cdots z_\g 
\nonumber \\
&=& y_\b \cdots y_2 z_2 \cdots z_\g 
\ + \ y_\b \cdots y_2 P_0^\perp(y_1z_1) z_2 \cdots z_\g.
\nonumber
\end{eqnarray}
Continuing in this way, we obtain, even for $\d = 1$ or
$\mu(\py) \neq \mu(\pz)$:

\setcounter{lemma}{1}
\begin{lemma}
\label{lem6.2}
Let $\py, \pz \in \cB$ with $\py = y_1 \cdots y_\b$
and $\pz = z_1 \cdots z_\g$, and let $\d \geq 1$ be
the integer such that $y_i = z_i$ for all $i < \d$
while $y_\d \neq z_\d$ (including the case 
$\b = \d-1$ or $\g = \d-1$). Then
$$
\py^*\pz = \sum_{i=1}^\d y_\b \cdots y_{i+1}
P_0^\perp(y_i z_i)z_{i+1} \cdots z_\g,
$$
where
\begin{enumerate}
\item One should replace $P_0^\perp(y_iz_i)$ by $1$ if
$\mu(y) \neq \mu(z)$ so that $\py^*\pz \in \cB$.
\item One should replace the summand for $i = \d$ by $1$
if $y_\d$ and $z_\d$ are both not present, i.e. if $\py=\pz$.
\item If $\b = \d - 1$ then no $y_j$'s should appear on
the left of the term for $i = \d$, 
and similarly if $\g = \d -1$.
\end{enumerate} 
\end{lemma}  

Suppose now that for some $i \leq \d$ we have
$$
P_m(y_\b \cdots y_{i+1}P_0^\perp(y_iz_i)z_{i+1} \cdots
z_\g) \neq 0.
$$
Then there must be an $r \in \cB^{\mu(y_i)}$ , $r \neq 1$,
such that $\s(ry_iz_i) \neq 0$ and
$$\ell(y_\b \cdots y_{i+1}) + \ell(r) + \ell(z_{i+1}
\cdots z_\g) = m.
$$
But because $\s(ry_iz_i) \neq 0$ we also have, by the
properties of filtrations,
\[
\ell(y_i) + \ell(z_i) \geq \ell(r) 
\geq |\ell(y_i) - \ell(z_i)|.
\]
Thus
\[
\ell(y_\b \cdots y_i) + \ell(z_i \cdots z_\g)
\geq \ell(y_\b \cdots y_{i+1}) + |\ell(y_i) - \ell(z_i)|
+ \ell(z_{i+1} \cdots z_\g).
\]
Let $\pw = y_1 \cdots y_{i-1} = z_1 \cdots z_{i-1}$. It 
follows from above that
\[
\min\{\ell(y_1 \cdots y_i), \ell(z_1 \cdots z_i)\}
\geq (\ell(\py) + \ell(\pz) -m)/2 \geq \ell(\pw).
\]
Recall that $q = (k+n-m)/2$. Since $\ell(\py) = k$
and $\ell(\pz) = n$, we see that 
$\ell(\pw) = \ell(y_1 \cdots y_{i-1}) \leq q$,
while $\ell(y_1 \cdots y_i) \geq q$ so that
\[\ell(y_\b \cdots y_{i+1}) \leq k-q.
\]
Similarly
\[
\ell(z_{i+1} \cdots z_\g) \leq n-q.\]
Notice that
\[
(k-q) + (n-q) = m,
\]
so that we can not have simultaneously 
$\ell(y_\b \cdots y_{i+1}) = k-q$ and 
$\ell(z_{i+1} \cdots z_\g) = n-q$. We 
summarize the above observations by:

\setcounter{lemma}{2}
\begin{lemma}
\label{lem6.3}
Suppose that $\py$ and $\pz$ are such that 
$\mu(\py) = \mu(\pz)$. let $\d$ be as defined as
above. If for some $i \leq \d$ we have
\[P_m(y_\b \cdots y_{i+1}P_0^\perp(y_iz_i)z_{i+1} 
\cdots z_\g) \neq 0,
\]
then $\py$ and $\pz$ are of the form $\py = \pw^* u \hat \ps$
and $\pz = \pw^* v \hat \pt$ where $\ell(\pw) \leq q$,
$\ell(\hat \ps) \leq k-q$, $\ell(\hat \pt) \leq n-q$, and
$u, v \in \cB^1\cup \cB^2$
with $\mu(\pw) \neq \mu(u) \neq \mu(\hat \ps)$ and
$\mu(\pw) \neq \mu(v) \neq \mu(\hat \pt)$.
At least one of $u,v$ is not $1$, and if $u = 1$ then
also $\hat \ps = 1$, and similarly for $v$. Specifically,
$\pw = y_{i-1} \cdots y_1 = z_{i-1} \cdots z_1$ and $u = y_i$
and $v= z_i$, while
$\hat \ps = y_{i+1} \cdots y_\b$ and 
$\hat \pt = z_{i+1} \cdots z_\g$.
If $\b \leq i$ then $\hat \ps = 1$, 
and similarly for $\g \leq 1$.
Then
\[
P_m(y_\b \cdots y_{i+1}P_0^\perp(y_iz_i)z_{i+1} \cdots z_\g)
= P_m(\hat \ps^*P_0^\perp(uv)\hat \pt).
\]
Either $\ell(\hat \ps) < k-q$ or $\ell(\hat \pt) < n-q$ 
(or both). 
\end{lemma}

In order to be in a position to apply our assumption
that $(A^1, \s ^1)$ and ($A^2, \s^2)$ satisfy a 
Haagerup-type condition, we need to consider collectively
all the $\px$'s which may occur in the support of a fixed
term $y_\b \cdots y_{i+1}P_0^\perp(y_iz_i)z_{i+1} 
\cdots z_\g$. For this purpose it is convenient to
assume now that both $k-q \neq 0$ and $n-q \neq 0$.
At the end of the proof we will give separately the argument
for the remaining cases.
We also need to divide the situation into two cases,
depending on the structure of the $\px$'s. Let 
$\px = x_1 \cdots x_\a$. For the first case we assume
that there is a $j$ such that 
$\ell(x_1 \cdots x_j) < k-q$ while 
$\ell(x_1 \cdots x_{j+1}) > k-q$. (This will always
happen if $q$ is not an integer.) Thus we can express
$\px$ as $\px = \ps^*r\pt$ where $\mu(\ps) = \mu(\pt)$
and $\mu(r) \neq \mu(\ps)$, with $\ell(\ps) <k-q$
and $\ell(\pt) < n-q$. The second case will be that in 
which there is a $j$ such that $\ell(x_1 \cdots x_j) = k-q$.

\setcounter{notation}{3}
\begin{notation}
\label{not6.4}

Assume that $k-q \neq 0$ and $n-q \neq 0$. For any 
pair $(\ps, \pt)$ of elements of $\cB$ such that 
$\ell(\ps) < k-q$ and $\ell(\pt) < n-q$ we set:
\begin{itemize}
\item[a)] If $\mu(\ps) = \mu(\pt)$ (with $\ps = 1$
and/or \ $\pt = 1$ permitted --- recall our convention
about $\mu(1)$), then
\[
\cBst = \{\px \in \cB_m: \px = \ps^*r\pt, 
r \in \cB^1\cup \cB^2 \setminus \{1\},
\textrm{ and } \mu(\ps) \neq \mu(r) \neq \mu(\pt)\}.
\]
We let $E_{\ps, \pt}$ denote the linear span of 
$\cBst$, and we let $\Pst$ denote
the projection onto $E_{\ps,\pt}$.
\item[b)] If $q$ is an integer and $\mu(\ps) \neq \mu(\pt)$ 
(with $\ps = 1$
and/or \ $\pt = 1$ permitted), then
\begin{eqnarray}
\cC_{\ps,\pt} &=& \{\px \in \cB_m:\px = \ps^* r_1r_2\pt,
\ r_1 \in \cB^{\nu(r_2)}_{k-q-\ell(\ps)}, 
\ r_2 \in \cB^{\nu(r_1)}_{n-q-\ell(\pt)}, 
\nonumber\\
& & \quad \quad \mu(r_1) \neq \mu(\ps), \
\mu(r_2) \neq \mu(\pt)\}.
\nonumber
\end{eqnarray}
(Note that $\ell(r_i) \geq 1$ for $i=1,2$ since
$\ell(\ps) <k-q$ and $\ell(\pt) < n-q$.) We let 
$F_{\ps,\pt}$ denote the linear span of $\cC_{\ps,\pt}$
and we let $\Qst$ denote the projection onto 
$F_{\ps,\pt}$.
\end{itemize}
\end{notation}

\setcounter{lemma}{4}
\begin{lemma}
\label{lem6.5}
$\cB_m$ is the disjoint union of all the $\cB_{\ps,\pt}$'s
and $\cC_{\ps,\pt}$'s.
\end{lemma}

\begin{proof}
It is evident that the $\cB_{\ps,\pt}$'s are disjoint
among themselves, as are the $\cC_{\ps,\pt}$'s. If
$\px \in \cB_{\ps,\pt}$ for some $(\ps,\pt)$ then $\px$ is
not of the form $\pu\pv$ where $\pu \in \cB_{k-q}$
and $\pv \in \cB_{n-q}$, whereas all elements of 
any $\cC_{\ps,\pt}$ are of this form. Thus the 
$\cBst$'s are disjoint from the $\cCst$'s.

Let $\px \in \cB_m$ with $\px = x_1 \cdots x_\a$. 
Recall our assumption that $m \geq 1$. If $\px$ satisfies the 
conditions for the first case discussed just 
before Notation \ref{not6.4}, then $\px \in \cBst$
for the choice of $\ps,\pt$ given there. Suppose instead
that $\px$ does not satisfy the conditions of the first
case. Then there is a $j$ such that 
$\ell(x_1 \cdots x_j) = k-q$. (Thus $q$ is an integer.)
Since $k \neq q$, $\ell(x_j) \geq 1$. Thus we can write 
$x_1 \cdots x_j = \ps^*r_1$ with $r_1 = x_j$,
so $\ell(r_1) \geq 1$ 
and $\ell(\ps) + \ell(r_1) = k-q$, and 
$r_1 \in \cB^{\nu(\ps)}$ unless $\ps = 1$. Since
$(k-q)+(n-q) = m$, we will also have 
$\ell(x_{j+1} \cdots x_\a) = n-q \neq 0$, so that 
$x_{j+1} \cdots x_\a = r_2\pt$ with $\ell(r_2) \geq 1$,
$\ell(r_2) + \ell(\pt) = n-q$ and $r_2 \in \cB^{\nu(r_1)}$,
and $r_2 \in \cB^{\nu(\pt)}$ unless $\pt = 1$. Thus
$\px \in \cCst$ for this choice of $(\ps,\pt)$.
\end{proof}

\setcounter{corollary}{5}
\begin{corollary}
\label{cor6.6}
Assume that $k\neq q$ and $n\neq q$. Then
\[
P_m = (\bigoplus_{\mu(\ps) = \mu(\pt)}\Pst)
\oplus(\bigoplus_{\mu(\ps) \neq \mu(\pt)}\Qst),
\]
where $\ps = 1$ and $\pt = 1$ are permitted.
\end{corollary}

As this corollary suggests, we will now examine
$\Pst(a^*\xi)$ and $\Qst(a^*\xi)$ in order to 
obtain the estimate we need for $P_m(a^*\xi)$.

\setcounter{lemma}{6}
\begin{lemma}
\label{lem6.7}
Let $(\ps,\pt)$ be such that $\mu(\ps) = \mu(\pt)$,
with $\ps = 1$ and $\pt = 1$ permitted. Let $\py \in \cB_k$
and $\pz \in \cB_n$ be given. If $\Pst(\py^*\pz) \neq 0$,
then $\py$ and $\pz$ are of the form $\py = \pw^* u \ps$
and $\pz = \pw^* v \pt$ where 
\begin {itemize}
\item[] $u, v \in \cB^1\cup\cB^2 \setminus \{1\}$ and
$\mu(u) = \mu(v)$,
\item[] $\mu(\ps) \neq \mu(u) \neq \mu(\pw)$ and 
$\mu(v) \neq \mu(\pt)$, 
\item[]$\ell(\pw) \leq q$, with $\pw = 1$ permitted. 
\end{itemize}
(Consequently $\ell(u) = k - \ell(\ps) - \ell(\pw)$ and
$\ell(v) = n - \ell(\pt) - \ell(\pw)$.) 
\newline Then
\[
\Pst(\py^*\pz) = \ps^*P_{m(\ps,\pt)}(uv)\pt,
\]
where $m(\ps,\pt) = m -\ell(\ps) - \ell(\pt)$.
\end{lemma}

\begin{proof}
This follows from Lemma \ref{lem6.3} when we set
$\ps = \hat \ps$ and $\pt = \hat \pt$ there.
\end{proof}

\setcounter{lemma}{7}
\begin{lemma}
\label{lem6.8}
Let $(\ps,\pt)$ be such that $\mu(\ps) \neq \mu(\pt)$,
with $\ps = 1$ and $\pt = 1$ permitted. Let $\py \in \cB_k$
and $\pz \in \cB_n$ be given. If $\Qst(\py^*\pz) \neq 0$,
then $\py$ and $\pz$ are in one and only one of the forms:
\begin{itemize}
\item[a)] $\py = \pw^* u\ps$ and $\pz = \pw^* vr_2\pt$ where 
\begin{itemize}
\item[] $u, v \in \cB^1\cup\cB^2 \setminus \{1\}$ and 
$\mu(u) = \mu(v)$, 
\item[] $\mu(\ps) \neq \mu(u) \neq \mu(\pw)$,
and $\mu(v) \neq \mu(r_2) \neq \mu(\pt)$, 
\item[] $\ell(r_2\pt) = n-q$ while $\ell(\pw) \leq q$.
\end{itemize}
Then 
\[
\Qst(\py^*\pz) = \ps^*P_{m(\ps,r_2\pt)}(uv)r_2\pt,
\]
where $m(\ps, r_2\pt) = m - \ell(\ps) - \ell(r_2\pt)$.
\item[b)] $\py = \pw^* ur_1\ps$ and $\pz = \pw^* v\pt$ 
with similar restrictions as above, and $\ell(r_1\ps) = k-q$
while $\ell(\pw) \leq q$. Then
\[
\Qst(\py^*\pz) = \ps^*r_1P_{m(\ps^*r_1,\pt)}(uv)\pt.
\]
\end{itemize}
\end{lemma}

\begin{proof}
In the notation of Lemma \ref{lem6.3} this is the case
in which either $\ell(\hat \ps) = k-q$ or 
$\ell(\hat \pt) = n-q$ (but not both). 
If $\ell(\hat \pt) = n-q$ with $\hat \pt = z_{i+1} \cdots z_\g$, 
then we must have $r_2 = z_{i+1}$ and 
$\pt = z_{i+2} \cdots z_\g$. We also have $\ps = \hat \ps$.
This gives case a). If, instead, $\ell(\hat \ps) = n-q$,
then we must have $r_1 = y_{i+1}$ and 
$\ps = y_{i+2} \cdots z_\g$, while $\pt = \hat \pt$. This
gives case b).
\end{proof}

\begin{proof}
Proof of theorem \ref{thm6.1}. Suppose now that $a \in E_k$
and $\xi \in E_n$, with $a = \sum a(\py)\py$ and
$\xi = \sum \xi(\pz)\pz$. Let $(\ps,\pt)$ be such that
$\Pst$ is defined. For any $\ps'$ and $\pw \in \cB$ set
$k(\ps', \pw) = k -\ell(\ps') - \ell(\pw)$ and
$n(\ps', \pw) = n - \ell(\ps') - \ell(\pw)$. Then from
Lemma \ref{lem6.7} we have
\[
\Pst(a^*\xi)  =  \sum_{\pw, u, v}
\bar a(\pw^*u\ps)\xi(\pw^*v\pt)\ps^*P_{m(\ps,\pt)}(uv)\pt
\]
where in the above sum
\begin{itemize}
\item[] $u, v \in \cB^1\cup\cB^2 \setminus \{1\}$ and
$\mu(u) = \mu(v)$,
\item[] $\mu(\ps) \neq \mu(u)\neq \mu(\pw)$ and 
$\mu(v) \neq \mu(\pt)$, 
\item[] $\ell(\pw) \leq q, \ \ell(u) = k(\ps,\pw)$, and 
$\ell(v) = n(\pt,\pw)$.
\end{itemize}
This sum can be rewritten as
\[\ps^*\big(\sum_{\substack{\ell(\pw) \leq q \\
\mu(\pw) = \mu(\ps)}} P_{m(\ps,\pt)}(\tilde a_{\ps,\pw} \
\tilde \xi_{\pt,\pw})\big)\pt,
\]
where we have set
\begin{eqnarray}
\tilde a_{\ps,\pw} & = 
& \sum_{u \in \cB^{\nu(\ps)}_{k(\ps,\pw)}} \bar a(\pw^*u\ps)u 
\nonumber \\
\tilde \xi_{\pt,\pw} & = 
& \sum_{v \in \cB^{\nu(\pt)}_{n(\pt,\pw)}} \xi(\pw^*v\pt)v. 
\nonumber
\end{eqnarray}
Note that for any $\px \in \cB$ and $b \in A^{\nu(\px)}$
we have $\|b\px\|_2 = \|b\|_2 = \|\px^*b\|_2$.
Consequently
\begin{eqnarray}
\|\Pst(a^*\xi)\|^2_2 & \leq & 
\big(\sum_{\pw}\|P^{\nu(\ps)}_{m(\ps,\pt)}(\tilde a_{\ps,\pw} \
\tilde \xi_{\pt,\pw)})\|_2\big)^2 \nonumber \\
& \leq & \big(\sum_{\pw} C\|\tilde a_{\ps,\pw}\|_2\|
\tilde \xi_{\pt,\pw}\|_2\big)^2 \nonumber \\
& \leq & C^2\big(\sum_{\pw}\|\tilde a_{\ps,\pw}\|^2_2\big)
\big(\sum_{\pw'}\|\tilde \xi_{\pt,\pw'}\|^2_2\big)
\nonumber \\
& = &C^2\big(\sum_{\pw}\sum_u |a(\pw^*u\ps)|^2\big)
\big(\sum_{\pw'}\sum_v |\xi(\pw'^*v\pt|^2\big).
\nonumber
\end{eqnarray}

The second inequality is the crucial place where we use
the assumption that $A^1$ and $A^2$ satisfy a
Haagerup-type condition with constant $C$. The third
inequality comes from the Cauchy-Schwarz inequality. 

We have seen that the $\Pst$'s form an orthogonal family
of projections. Consequently, with the understanding that
$\ell(\ps) < k-q$ , $\ell(\pt) < n-q$, and
$\mu(\ps) = \mu(\pt)$, with $\ps = 1$ and $\pt = 1$
permitted, we obtain
\begin{eqnarray}
\|\sum_{\ps,\pt}\Pst(a^*\xi)\|^2_2 & = &
\sum_{\ps,\pt}\|\Pst(a^*\xi)\|^2_2 \nonumber \\
& \leq & C^2\sum_{\ps,\pt}\big(\sum_{\pw}\sum_u |a(\pw^*u\ps)|^2\big)
\big(\sum_{\pw'}\sum_v |\xi(\pw'^*v\pt|^2\big).
\nonumber
\end{eqnarray}
Now any given $\py \in \cB_k$ has a unique expression as
$\py = \pw u \ps$ for some $\pw$ with $\ell(\pw) \leq q$
and some $\ps$ with $\ell(\ps) < k-q$, and similarly
for $\pz \in \cB_n$ as $\pz = \pw v \pt$. It is easily 
seen from this that we obtain
\[
 \|\sum_{\ps,\pt}\Pst(a^*\xi)\|^2_2 \leq C^2\|a\|^2_2
\|\xi\|^2_2.
\]
Notice that if $q$ is not an integer, so that
$P_m = \sum \Pst$, then this already gives the 
desired inequality, and the proof of the theorem is complete.

Suppose instead that $q$ is an integer and that 
$\mu(\ps) \neq \mu(\pt)$, so that $\Qst$ is defined.
Then from Lemma \ref{lem6.8} we have
\begin{eqnarray}
\Qst(a^*\xi) & = & \sum_{\pw, u, v, r_2}
\bar a(\pw^*u\ps)\xi(\pw^*vr_2\pt)\ps^*
P_{m(\ps,r_2\pt)}(uv)r_2\pt \nonumber \\
& + & \sum_{\pw, u, v, r_1}
\bar a(\pw^*ur_1\ps)\xi(\pw^*v\pt)\ps^*r_1
P_{m(r_1\ps,\pt)}(uv) \pt, \nonumber
\end{eqnarray}
where in both sums $\pw \in \cB$ with $\ell(\pw) \leq q$
and $u, v, r_1, r_2 \in \cB^1\cup\cB^2 \setminus \{1\}$
with $\mu(u) = \mu(v)$,
while in the first sum
\begin{itemize}
\item[] $\mu(\ps) \neq \mu(u) \neq \mu(\pw)$ and
\ $\mu(v) \neq \mu(r_2) \neq \mu(\pt)$, 
\item[] $\ell(u) = k(\ps, \pw), \ \ell(v) = n(r_2\pt, \pw)$,
and $\ell(r_2\pt) = n-q$,
\end{itemize}
whereas in the second sum
\begin{itemize}
\item[] $\mu(\ps) \neq \mu(r_1) \neq \mu(u)$ and
$\mu(\pw) \neq \mu(v) \neq \mu(\pt)$, 
\item[] $\ell(u) = k(r_1\ps, \pw), \
\ell(v) = n(\pt,\pw)$ and $\ell(r_1\ps) = k-q$.
\end{itemize}
For each $\pw$ with $\ell(\pw) \leq q$ let us define 
$\tilde a_{\ps,\pw}$, etc. much as before by
\begin {eqnarray}
\tilde a_{\ps,\pw} & = & \sum_u a(\pw^*u\ps)u, \quad
\tilde \xi_{r_2\pt,\pw} = \sum_v\xi(\pw^*vr_2\pt)v,
\nonumber \\
\tilde a_{r_1\ps,\pw} & = & \sum_u a(\pw^*ur_1\ps)u, \quad
\tilde \xi_{\pt,\pw} = \sum_v\xi(\pw^*v\pt)v,
\nonumber
\end{eqnarray}
with the restrictions on $u$ and $v$ as above. Then in 
terms of this notation we have
\[
\Qst(a^*\xi) = \sum_{\pw,r_2}\ps^*P_{m(\ps,r_2\pt)}
(\tilde a_{\ps,\pw}\tilde \xi_{r_2\pt,\pw})r_2\pt
+ \sum_{\pw,r_1}\ps^*r_1P_{m(r_1\ps,\pt)}
(\tilde a_{r_1\ps, \pw}\tilde \xi_{\pt,\pw})\pt.
\]
Since the two summands above may not be orthogonal, but
the terms within each sum over $r_1$ and $r_2$ are 
orthogonal, we obtain
\begin{eqnarray}
\|\Qst(a^*\xi)\|^2_2 & \leq & 2\big(\|\sum_{\pw,r_2}
\ps^*P_{m(\ps,r_2\pt)}(\tilde a_{\ps,\pw}
\tilde \xi_{r_2\pt,\pw})r_2\pt\|^2_2 
\nonumber \\ 
&{}&  \quad\quad  +\ \
\|\sum_{\pw,r_1}\ps^*r_1P_{m(r_1\ps,\pt)}
(\tilde a_{r_1\ps,\pw}\tilde \xi_{\pt,\pw})\pt
\|^2_2\big)
\nonumber \\
& \leq & 2\sum_{r_2}\big(\sum_{\pw}
\|P_{m(\ps,r_2\pt)}(\tilde a_{\ps,\pw}
\tilde \xi_{r_2\pt,\pw}) \|_2\big)^2 
\nonumber \\ 
&{}& \quad \quad + 
2\sum_{r_1}\big(\sum_{\pw}\|P_{m(r_1\ps,\pt)}
(\tilde a_{r_1\ps,\pw}\tilde \xi_{\pt,\pw})\|_2\big)^2
\nonumber \\
& \leq & 2\sum_{r_2}\big(\sum_{\pw} C\|\tilde a_{\ps,\pw}\|_2
\|\tilde \xi_{r_2\pt,\pw}\|_2\big)^2
\nonumber \\
&{}& \quad \quad +\ \ 2\sum_{r_1}\big(\sum_{\pw} 
C\|\tilde a_{r_1\ps,\pw}
\|_2\|\tilde \xi_{\pt,\pw}\|_2\big)^2
\nonumber \\
& \leq & 2C^2\big(\sum_{\pw}
\|\tilde a_{\ps,\pw}\|^2_2\big)
\big(\sum_{\pw',r_2}\|\tilde \xi_{r_2\pt,\pw'}\|^2_2\big)
\nonumber \\
&{}& \quad \quad +\ \ 2C^2\big(\sum_{\pw, r_1}
\|\tilde a_{r_1\ps,\pw}\|^2_2\big)
\big(\sum_{\pw'}\|\tilde \xi_{\pt,\pw'}\|^2_2\big) .
\nonumber
\end{eqnarray}
We have seen that the $\Qst$'s form an orthogonal family 
of projections. Consequently, with the understanding
that $\ell(\ps) < k-q$, \ $\ell(\pt)<n-q$ and 
$\mu(\ps) \neq \mu(\pt)$, with $\ps = 1$ and/or
$\pt = 1$ permitted, we obtain
\begin{eqnarray}
\|\sum_{\ps,\pt}\Qst(a^*\xi)\|^2_2 
&=& \sum_{\ps,\pt}\|\Qst(a^*\xi)\|^2_2
\nonumber \\
& \leq & 2C^2\sum_{\ps,\pt}
\Big((\sum_{\pw}\|\tilde a_{\ps,\pw}\|^2_2)
(\sum_{\pw',r_2}\|\tilde \xi_{r_2\pt,\pw'}\|^2_2)
\nonumber \\
&{}&\quad \quad +\ \ (\sum_{\pw,r_1}
\|\tilde a_{r_1\ps,\pw}\|^2_2)
(\sum_{\pw'}\|\tilde \xi_{\pt,\pw'}\|^2_2)\Big).
\nonumber
\end{eqnarray}
Now again any given $\py \in \cB_k$ has a unique expression
as $\py = \pw u\ps$ for some $\pw$ with 
$\ell(\pw) \leq q$ and some $\ps$ with 
$\ell(\ps) < k-q$; and furthermore, if $\py$ can be
expressed as $\py = \pw ur_1\ps$ with 
$\ell(r_1) = k-1-\ell(\ps)$ and $\ell(u)+\ell(\pw)=q$,
then this expression too is unique. A similar statement 
holds for any $\pz \in \cB_n$ as $\pz = \pw v\pt$ or
$\pz = \pw vr_2\pt$. In the same way as for the $\Pst$'s
it is then easily seen that
\[
\|\sum_{\ps,\pt}\Qst(a^*\xi)\|^2_2 \leq 
4C^2\|a\|^2_2\|\xi\|^2_2.
\]
Since $P_m$ is the orthogonal sum of the $\Pst$'s and
the $\Qst$'s, it follows that
\[
\|P_m(a^*\xi)\|_2 \leq \sqrt{5}C\|a\|_2\|\xi\|_2
\]
as desired.

Finally, we must treat the cases in which $k-q=0$ or $n-q=0$.
If $k-q=0$ \ then $m+k=n$. We follow the pattern of proof of
the previous cases, and so allow ourselves less detailed
notation and discussion. For any $\pt \in \cB$ with 
$\ell(\pt) < m$ set
\[
\cB_{(\pt)} = \{\px\in \cB_m:\px = r\pt \textrm{ with }
r \in \cB^1\cup\cB^2 \setminus \{1\}, \ \mu(r)\neq \mu(\pt)\}.
\]
We permit $\pt = 1$. It is easily seen that the 
$\cB_{(\pt)}$'s are disjoint and that their union is
$\cB_m$. We let $E_{(\pt)}$ denote the linear span of
$\cB_{(\pt)}$, and we let $P_{(\pt)}$ denote the projection
onto $E_{(\pt)}$.

\setcounter{lemma}{8}
\begin{lemma}
\label{lem6.9}
Let $\py \in \cB_k$ and $\pz \in \cB_n$. If
$P_{(\pt)}(\py^*\pz) \neq 0$ then $\py$ and $\pz$ are
of the form $\py = \pw^*u$ and $\pz = \pw^*v\pt$ where
$u, v \in \cB^1\cup \cB^2$, with $\ell(\pw) \leq k$ and
$\mu(u) \neq \mu(\pw) \neq 
\mu(v) \neq \mu(\pt)$ and $v \neq 1$. 
(But we may have $u=1$.) Then $P_{(\pt)}(\py^*\pz) = 
P_{m(\pt)}(uv)\pt$ where  $m(\pt) = m - \ell(\pt)$.
\end{lemma}

\begin{proof}
According to Lemma \ref{lem6.3} we can express $\py$ 
and $\pz$ as $\py = \pw^*u\hat \ps$ and 
$\pz = \pw^*v\hat \pt$ where among the conditions we
have $\ell(\hat \ps) \leq k-q = 0$. Thus $\hat \ps = 1$. So
$\py = \pw^*u$ with $\mu(\pw) \neq \mu(u)$. We will also
have $\ell(\pw) \leq q = k$ and $\ell(\hat \pt) \leq n-q=m$.
Suppose that $v=1$. Then 
$\ell(\pw) + \ell(\hat \pt) = \ell(\pz) = n = k+m$, and
so $\ell(\pw) = k$, $\ell(\hat \pt) = m$ and $u=1$,
which contradicts Lemma \ref{lem6.3}. Thus $v \neq 1$. We can
set $\pt = \hat \pt$. Then from Lemma \ref{lem6.3} we have
$P_{(\pt)}(\py^*\pz) = 
P_{m(\pt)}(uv)\pt$.
\end{proof}

Suppose now that $a \in E_k$ and $\xi \in E_n$. 
Then, much as in 
the previous cases, we have
\[
P_{(\pt)}(a^*\xi) = 
\sum_{\pw, u, v}\bar a(\pw^*u)\xi(\pw^*v\pt)
P_{m(\pt)}(uv)\pt,
\]
where the conditions on $\pw, u, v$ are as above. We set
\[
\tilde a_{\pw} = \sum_u\bar a(\pw^*u)u, \qquad
\tilde \xi_{\pt,\pw} = \sum_v\xi(\pw^* v \pt)v\pt.
\]
Thus
\begin{eqnarray}
\|P_{(\pt)}(a^*\xi)\|^2_2 &\leq& 
\big(\sum_{\pw}\|P_{m(\pt)}(\tilde a_{\pw}
\tilde \xi_{\pt,\pw})\|_2\big)^2
\nonumber \\
&\leq& \big(\sum_{\pw}C\|\tilde a_{\pw}\|_2
\|\tilde \xi_{\pt,\pw}\|_2\big)^2
\leq C^2\big(\sum_{\pw}\|\tilde a_{\pw}\|^2_2\big)
\big(\sum_{\pw'}\|\tilde \xi_{\pt,\pw'}\|^2_2\big)
\nonumber \\
&=& C^2\big(\sum_{\pw,u}|\bar a(\pw^*u|^2\big)
\big(\sum_{\pw',v}|\xi(\pw'^*v\pt|^2\big).
\nonumber
\end{eqnarray}
Consequently
\begin{eqnarray}
|P_m(a^*\xi)|^2_2 &=& \sum_{\pt}\|P_{(\pt)}(a^*\xi)\|^2
\nonumber \\
&\leq&  C^2\sum_{\pt}
\big(\sum_{\pw,u}|\bar a(\pw^*u|^2\big)
\big(\sum_{\pw',v}|\xi(\pw'^*v\pt|^2\big).
\nonumber
\end{eqnarray}

Now because $k+m=n$ it is easily seen that any 
given $\pz \in \cB_n$ has a unique expression
as $\pz = \pw v\pt$ where $\ell(\pw) \leq k$, 
$\ell(\pt) < m$, and 
$v \in \cB^1\cup\cB^2 \setminus \{1\}$. However
a \ $\py \in \cB_k$  will have two expressions as
$\pw u$ with $\ell(\pw) \leq k$ and $u \in \cB^1\cup\cB^2$
(and $\mu(\pw) \neq \mu(u)$), one of which will 
be $\py = \pw$. It follows that
\[
\|P_m(a^*\xi)\|^2_2 \leq 2C^2\|a\|^2_2\|\xi\|^2_2,
\]
which implies the desired inequality.

Finally, we must deal with the case in which $n-q = 0$.
But this case follows from essentially the mirror image
of the above argument, in which now for $\ell(\ps) <m$ the 
elements of $\cB_{(\ps)}$ have form $\px = \ps r$, and 
later we find that $(\py,\pz)$ must have the form
$\py = \pw^*u\ps$ and $\pz = \pw^*v$.
\end{proof}

\setcounter{question}{9}
\begin{question}
\label{que6.10}
What happens for amalgamated free products of $C^*$-algebras?
What happens if $\s_1$ and $\s_2$ are not tracial?
\end{question}


\end{document}